# Uncertainty modeling method for wind and solar power output in building integrated energy systems under continuous anomalous weather

Deyi Shao[a], Hongru Li[a, *], Jingsheng Li[a], Xia Yu[a], Xiaoyu Sun[a], Bowen Han[a]

[a] *College of Information Science and Engineering, Northeastern University, NO. 3-11 Wenhua Road, Shenyang, 110819, China.*

**Abstract**

The increasing occurrence of continuous anomalous weather events has intensified the uncertainty in wind and photovoltaic power generation, posing significant challenges to the operation and optimization of building integrated energy systems. Existing studies often neglect the interdependence between successive anomalous weather events and their collective impact on wind and solar power output. Additionally, conventional modeling approaches struggle to accurately capture the nonlinear fluctuations induced by these weather conditions. To address this gap, this study proposes an uncertainty modeling method based on stochastic optimization and scenario generation. The Weibull and Beta distributions characterize the probabilistic properties of wind speed and solar irradiance, respectively, while the Copula function captures the dependence between wind speed and precipitation, enabling the construction of a wind-solar power uncertainty model that incorporates the joint distribution of consecutive anomalous weather events. A Monte Carlo-based scenario generation approach is employed to construct a dataset representing anomalous weather characteristics, followed by a probabilistic distance-based scenario reduction technique to enhance modeling efficiency. Furthermore, the unscented transformation method is introduced to mitigate nonlinear propagation errors in wind and solar power state estimation. Case studies demonstrate that the proposed method effectively characterizes the fluctuation patterns of wind and solar power under continuous anomalous weather conditions while preserving the statistical properties of the original data. These findings provide a reliable basis for improving the operational resilience of building integrated energy systems under extreme weather scenarios.

## 1. Introduction

Against the backdrop of a global transition toward a low-carbon energy structure, the share of renewable energy, particularly wind and photovoltaic power, in building integrated energy systems (BIES) continues to rise. However, intensifying global climate change has led to a significant increase in both the frequency and severity of anomalous climate events. Given the strong dependence of wind and solar resources on weather conditions, their power outputs exhibit high randomness and volatility. Anomalous climate events encompass both extreme weather and continuous anomalous weather [1]. The former refers to rare and severe weather

* Corresponding author.
*E-mail address:* shaodeyi@stumail.neu.edu.cn (D. Shao), lihongru@ise.neu.edu.cn (H. Li), 2300893@stu.neu.edu.cn (J. Li), yuxia@ise.neu.edu.cn (X. Yu), sunxiaoyu1@ise.neu.edu.cn (X. Sun), 2200821@stu.neu.edu.cn (B. Han).

| | |
|---|---|
| **Nomenclature** | |
| AIC | Akaike Information Criterion |
| GEV | Generalized Extreme Value |
| GP | Generalized Pareto |
| BIES | Building Integrated Energy Systems |
| K-S | Kolmogorov-Smirnov |
| PV | Photovoltaic system |
| QQ | Quantile-Quantile |
| UT | Unscented Transform |
| WT | Wind turbine |

conditions, while the latter describes sustained deviations of weather states from their mean conditions, potentially involving the coupling and superposition of multiple anomalous weather elements. The complexity and uncertainty of anomalous climate events, particularly their persistence and extremity, pose severe risks to human activities and the stable operation of BIES [2]. Therefore, when planning BIES under continuous anomalous weather conditions, considering the uncertainty of wind-solar power output is crucial to ensuring system stability and reliability [3]. In this context, investigating the nonlinear fluctuations of wind-solar power output induced by anomalous weather, analyzing the interdependencies among anomalous weather events, and assessing their comprehensive impact on wind-solar power output are of great significance [4].

*1.1. Literature review*

The nonlinear fluctuations in wind-solar power output induced by anomalous weather pose significant challenges for uncertainty modeling. To address this issue, stochastic programming, scenario analysis, robust optimization, and interval analysis have been widely applied in wind-solar power output uncertainty modeling [5]. Most researchers focus on accurately characterizing the fluctuation characteristics of wind-solar power output to enhance the stability and economic efficiency of BIES under extreme meteorological conditions. Niu et al. proposed an annual system operation scenario generation method considering extreme meteorological events, improving the accuracy of renewable energy output scenarios and verifying their impact on system operation [6]. Li et al. developed a wind and solar power output uncertainty modeling method based on historical data, replacing the expected output with the upper α-quantile and integrating Monte Carlo simulation to optimize power system planning, effectively reducing load loss rates and cost risks [7]. Perera et al. combined robust optimization with multi-climate scenario analysis to quantify the impact of extreme weather on renewable energy potential and power supply reliability, thereby enhancing the resilience of BIES against wind-solar power output uncertainty [8]. Han et al. proposed a wind turbine power curve modeling method based on interval extreme value probability density, which eliminates the need for abnormal data cleaning, enhances modeling accuracy and efficiency, and optimizes wind power forecasting [9]. They further leveraged duality theory and second-order cone programming to improve the reliability of load restoration. However, stochastic programming relies on a large number of sample simulations of random variables, often yielding only upper and lower bounds of variables. The accuracy of uncertainty modeling using scenario analysis depends on the number of selected scenarios and is limited to providing only the expected values of output variables [10].

Robust optimization considers multiple hazards by superimposing extreme events, while interval analysis assumes that the uncertainty parameters of power fluctuations caused by anomalous weather lie within known intervals. However, all these methods overlook the interdependencies among disasters and fail to establish correlation models [11]. In contrast, integrating stochastic programming with scenario analysis enables flexible handling of different scenarios, accurately reflecting the impact of anomalous weather and enhancing decision reliability and system resilience. Shahbazi et al. proposed a method combining scenario-based stochastic programming with bounded uncertainty robust optimization to simulate uncertainties in network equipment loading, energy prices, and availability under extreme weather conditions, and developed resilient building strategies to enhance distribution system robustness [12]. Despite these efforts to model wind-solar power output uncertainty under anomalous weather, they fail to address the nonlinear function propagation in power output state estimation. Anomalous weather involves multiple uncertain parameters, such as wind speed, solar radiation, temperature, and precipitation, necessitating an analysis of the impact of multiple random variables on wind-solar power output during uncertainty modeling. The unscented transform (UT) is an effective method for quantifying the uncertainty propagation of multiple random variables [13]. Athari et al. proposed a novel stochastic cascading failure model that employs the UT to analyze uncertainty parameters in power grid failures and utilizes an improved quasi-steady-state model to simulate voltage-dependent failures, validating its effectiveness under high wind power penetration [14]. Zhang et al. introduced an optimal management method based on flexible energy hubs, using the UT to model various uncertainties in loads, thereby reducing computational complexity and rapidly obtaining optimal solutions [15].

Despite the progress made in existing studies, most methods treat individual weather factors as independent variables, neglecting the interdependence among different weather events and their combined impact on wind-solar power output. In the generation of anomalous weather scenarios, there is a lack of systematic analysis of correlations between different weather events, resulting in limited accuracy in fitting the distribution of consecutive anomalous weather events. To capture the uncertainty of wind-solar power output under consecutive anomalous weather conditions, it is essential to incorporate the correlation analysis of different anomalous weather events into uncertainty modeling. The most accurate approach to reflect the dependency structure among anomalous weather events is to construct a joint distribution function that preserves the complete information from historical data [16]. Copula theory has been widely applied in constructing joint distribution models for multiple random variables [17]. Deng et al. introduced a spatiotemporal tail dependence structure and proposed a wind speed scenario generation method based on the C-vine Copula, effectively improving wind speed modeling accuracy and scenario generation reliability [18]. Sun et al. proposed a probabilistic solar power forecasting framework based on correlated weather scenario generation, utilizing Copula to model the joint distribution of weather variables, which significantly enhanced forecasting accuracy [19]. Yin et al. developed a joint distribution model for extreme persistent rainfall events using Copula, improving the accuracy of extreme rainfall identification and the calculation of design rainfall amounts [20]. Yang et al. employed a trivariate Copula method to assess the joint probability of storm surges and heavy rainfall induced by typhoons [21]. Additionally, parameter optimization was conducted using the particle swarm optimization algorithm to enhance the accuracy of disaster management. Chen et al. proposed a two-stage planning framework that captures correlations among extreme weather events via a Copula model, improving system flexibility and resilience while reducing system costs and outage rates [22]. Fu et al. introduced a Copula-based joint probability distribution method for wind speed and rainfall intensity, analyzing the impact of wind-rain combined loads on the failure probability of transmission lines [23]. Kazemi-Robati et al. proposed a scenario generation method based on Copula theory to account for the dependence structure among different random variables, addressing uncertainties in renewable energy generation and energy prices [24].

*1.2. Motivation and contributions*

The frequent occurrence of consecutive anomalous weather events exacerbates the uncertainty of wind-solar power output, posing significant challenges to the stability and reliability of BIES. As the primary sources of renewable energy, wind and photovoltaic power generation are highly influenced by meteorological conditions, and anomalous weather events further amplify this uncertainty, thereby affecting the scheduling and optimization efficiency of BIES. Therefore, accurately modeling the uncertainty of wind-solar power output, particularly under anomalous weather conditions, is crucial for enhancing the optimization and scheduling efficiency of BIES and strengthening their resilience to extreme climate challenges. By considering the complex fluctuations and multi-source uncertainties in wind-solar power output under consecutive anomalous weather conditions, more precise and comprehensive predictions can be provided for system scheduling and decision-making. Although existing studies have made progress in modeling wind-solar power output, most focus on individual meteorological factors without thoroughly exploring the intricate correlations among anomalous weather events. Furthermore, they have yet to effectively address the challenges of modeling the nonlinear fluctuations of wind-solar power output under anomalous weather conditions.

This study proposes a wind-solar power output uncertainty modeling method based on stochastic programming and scenario analysis. The aim is to systematically analyze the interdependencies and complex correlations among anomalous weather events, accurately capturing the nonlinear fluctuations of wind-solar power output and providing reliable data support for optimizing BIES scheduling under extreme climate conditions. The main contributions of this paper are as follows:

(1) Uncertainty modeling of wind-solar power output for BIES based on multiple stochastic variables: To address the strong randomness and nonlinear fluctuation characteristics of wind-solar power output under consecutive anomalous weather conditions, this study proposes a multi-stochastic-variable-based modeling approach. By integrating stochastic programming and the UT method, this approach characterizes the statistical properties of wind speed and solar irradiance using the Weibull and Beta distributions, respectively. The UT method is further employed to handle the nonlinear propagation of uncertainties, thereby enhancing the adaptability of the model to complex meteorological conditions.

(2) Joint distribution model of anomalous weather based on scenario analysis for BIES: A scenario tree of consecutive anomalous weather events is established based on scenario analysis. The Copula function is utilized to construct the joint distribution model of wind speed and precipitation, enabling an in-depth analysis of the mutual impact of consecutive anomalous weather events such as strong winds and heavy rainfall on wind-solar power output. By quantifying the correlation among anomalous weather events, the overall occurrence probability of consecutive anomalous weather is derived, providing a novel approach for evaluating the uncertainty in the scale of such weather events.

(3) Uncertainty quantification of wind-solar power output under consecutive anomalous weather conditions: A scenario generation method based on the joint distribution of anomalous weather is proposed. By integrating the Monte Carlo method, scenarios of wind-solar power output under anomalous weather conditions are generated. Additionally, a scenario reduction technique based on probability distance is applied to optimize the scenario set, ensuring that the generated scenarios not only preserve the critical features of anomalous weather impacts but also reduce computational complexity and mitigate the curse of dimensionality. Scenario evaluation is conducted to verify the representativeness of the generated scenarios, ensuring the provision of reliable data for wind-solar power output uncertainty modeling.

## 2. Methodology

Global climate change has intensified the frequency and intensity of anomalous weather events, leading to increased uncertainty in wind and solar output and causing significant fluctuations in energy demand, which poses challenges to the stable operation of BIES. At the same time, the operation of BIES, particularly the greenhouse gas emissions from fossil fuel consumption, further drives global warming and exacerbates the occurrence of anomalous weather events, creating a vicious cycle, as shown in Fig. 1. Anomalous weather not only affects key meteorological variables such as wind speed and solar irradiance, causing fluctuations in renewable energy output, but also increases the pressure on energy storage scheduling, thereby raising the risk to system security. Moreover, the imbalance between energy supply and demand under extreme climate conditions may exacerbate dependence on fossil fuels, further reinforcing the trend of climate warming. Therefore, there is an urgent need to develop accurate wind and solar output uncertainty modeling methods to quantify the impact of anomalous weather on renewable energy, providing theoretical support for the optimal scheduling and low-carbon transformation of BIES.

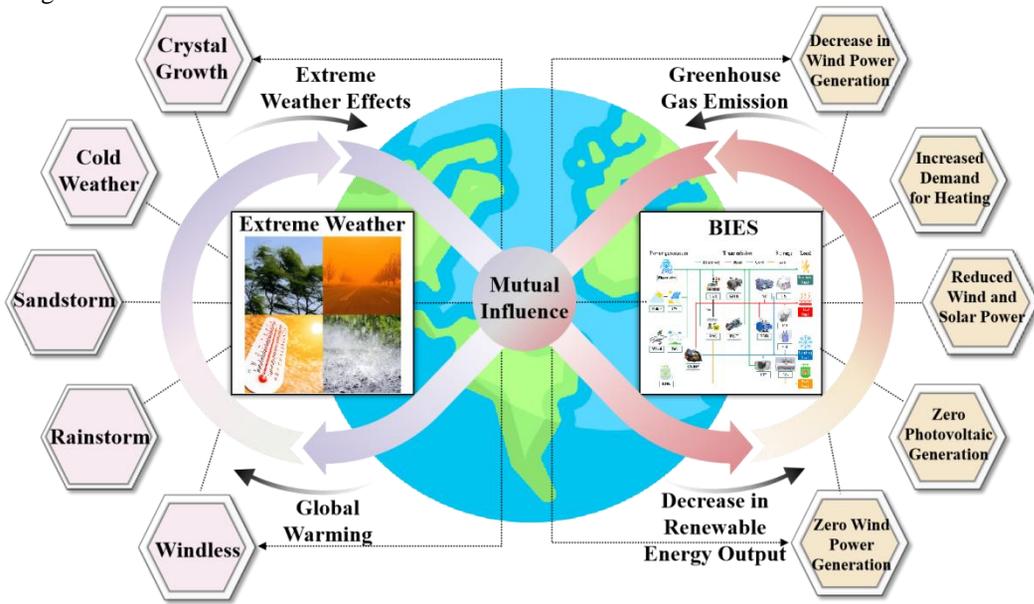

Fig. 1. The relationship between anomalous weather and BIES.

Fig. 2 illustrates the schematic of the scenario-based modeling method for wind and solar output uncertainty. The method comprehensively considers the nonlinear effects of anomalous weather events on wind speed and solar irradiance, using the Weibull and Beta distributions to depict the statistical characteristics of wind speed and solar irradiance, respectively. It also employs the Copula function to establish a joint distribution model for wind speed and precipitation, accurately describing the correlation of continuous anomalous weather events. Based on this, a scenario analysis approach is used to classify continuous anomalous weather events, and the Monte Carlo method is applied to generate a set of wind and solar output scenarios that account for the effects of anomalous weather. Additionally, probability-distance-based scenario reduction techniques are used to optimize the scenario set size, reducing computational complexity while maintaining the representativeness of the scenarios. Furthermore, to address the nonlinear propagation errors in wind and solar output, the Unscented Transform method is introduced to enhance the model's adaptability to complex meteorological conditions.

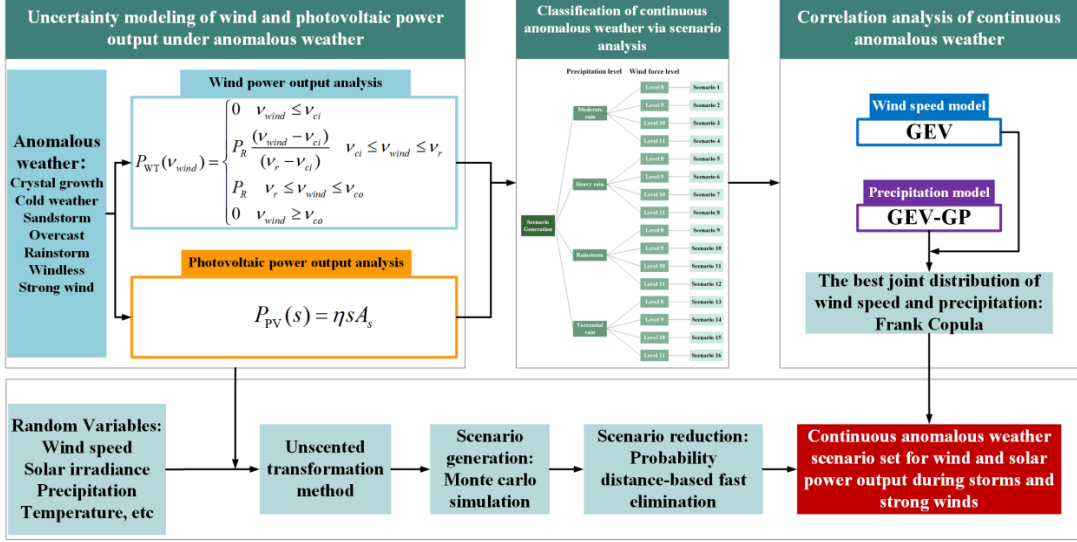

Fig. 2. Schematic of the scenario-based modeling method for wind and solar output uncertainty.

## 2.1. Wind and solar output uncertainty modeling

### 2.1.1. Wind turbine mode

In practical wind farms, the output of wind turbines (WT) depends on factors such as wind speed fluctuations, wind direction changes, equipment operating conditions, the wind energy utilization coefficient of the turbines, and the wake effects between turbines. The power generation of wind turbines is closely related to wind speed, constrained by both the magnitude of the wind speed and the randomness of wind speed changes. Meteorological data is suitable for assessing the wind energy capacity of a region. In the study of the probability distribution of wind speed forecasts, the two-parameter Weibull distribution is most widely used due to its reasonable performance. The probability density function and cumulative distribution of the Weibull distribution can be defined according to the following formulas [25]:

$$F(v_{wind}) = 1 - \exp\left[-\left(\frac{v_{wind}}{c}\right)^k\right] \tag{1}$$

$$f(v_{wind}) = \frac{k}{c}\left(\frac{v_{wind}}{c}\right)^{k-1} \exp\left[-\left(\frac{v_{wind}}{c}\right)^k\right] \tag{2}$$

Where $c$ is the scale parameter, reflecting the skewness of the curve; $k$ is the shape parameter, representing the predicted average wind speed of the region. When $k = 2$, the distribution becomes the standard Rayleigh distribution.

The shape parameter $k$ and scale parameter $c$ of wind speed can be calculated using Eqs. (3) and (4):

$$k = \left(\frac{v_w}{\mu_w}\right)^{-1.086} \quad (3)$$

$$c = \frac{\mu_w}{\Gamma(1+\frac{1}{k})} \quad (4)$$

Where $v_w$ and $\mu_w$ are the mean and standard deviation of wind speed, respectively; $\Gamma(\cdot)$ is the Gamma function.

In this study, taking the shape parameter $k = 2$, the relationship between the mean wind speed $v_w$ and the scale parameter $c$ is given by:

$$v_w = c\Gamma\left(1+\frac{1}{2}\right) = \frac{1}{2}c\Gamma\left(\frac{1}{2}\right) = \frac{\sqrt{\pi}}{2}c, \; c = \frac{2}{\sqrt{\pi}}v_w \quad (5)$$

Therefore, by substituting values into the probability density function, the cumulative distribution of the wind turbine model can be expressed as a function of the mean wind speed using Eqs. (1) and (2) (Kim et al., 2021):

$$F(V_{wind}) = 1 - \exp\left(-\left(\frac{\pi}{4}\right)\left(\frac{V_{wind}}{v_w}\right)^2\right) \quad (6)$$

$$f(V_{wind}) = \frac{\pi}{2}\frac{V_{wind}}{v_w^2}\exp\left(-\left(\frac{\pi}{4}\right)\left(\frac{V_{wind}}{v_w^2}\right)^2\right) \quad (7)$$

The output of the wind turbine is influenced by wind speed variations and exhibits cut-in and cut-out characteristics. When the wind speed is below the cut-in wind speed, the turbine is in a shutdown state and does not generate power. When the wind speed is between the cut-in and rated wind speeds, the turbine starts operating. For simplification, the relationship between wind power and wind speed in this phase is typically linearized. When the wind speed reaches the range between the rated wind speed and the cut-out wind speed, the turbine operates at rated power. If the wind speed exceeds the cut-out wind speed, the turbine will be forced to shut down to prevent structural damage.

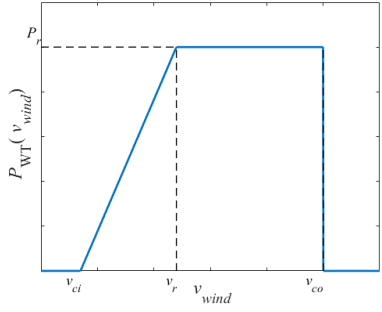

Fig. 3. Wind power versus wind speed curve.

For a given wind turbine, the wind power-speed curve shown in Fig. 3, the output power characteristics can be calculated using Eqs. (8) [26].

$$P_{WT}(v_{wind}) = \begin{cases} 0 & v_{wind} \leq v_{ci} \\ P_R \dfrac{(v_{wind} - v_{ci})}{(v_r - v_{ci})} & v_{ci} \leq v_{wind} \leq v_r \\ P_R & v_r \leq v_{wind} \leq v_{co} \\ 0 & v_{wind} \geq v_{co} \end{cases} \quad (8)$$

Where $v_{ci}$, $v_r$, $v_{co}$, and $v_{wind}$ represent the cut-in wind speed, rated wind speed, cut-out wind speed, and actual wind speed of the turbine, respectively, and PR represents the rated power of the turbine.

### 2.1.2. Photovoltaic power generation system model

The output of a photovoltaic power generation system is entirely dependent on solar irradiance, which is influenced by factors such as weather, location, day-night cycles, and seasons. The randomness of the solar irradiance distribution must be considered. Currently, the Beta distribution is commonly used for modeling its parameters [27], as shown in the following equations:

$$f(s) = \begin{cases} \dfrac{\Gamma(\alpha+\beta)}{\Gamma(\alpha)\Gamma(\beta)} s^{\alpha-1}(1-s)^{\beta-1} & 0 \leq s \leq 1, \alpha \geq 0, \beta \geq 0 \\ 0 & \text{otherwise} \end{cases} \quad (9)$$

$$f(s) = \int_0^s \dfrac{\Gamma(\alpha+\beta)}{\Gamma(\alpha)\Gamma(\beta)} s^{\alpha-1}(1-s)^{\beta} - 1 ds \quad (10)$$

Where $s$ represents the solar radiation (kW/m$^2$), and $\alpha$ and $\beta$ are the location and shape parameters of the Beta distribution, respectively. These parameters can be determined based on the average value and standard deviation of the solar irradiance, as shown in the following equation [28]:

$$\alpha = \left(\frac{\mu(1+\mu)}{\sigma^2} - 1\right) \tag{11}$$

$$\beta = (1-\mu)\left(\frac{\mu(1+\mu)}{\sigma^2} - 1\right) \tag{12}$$

Where $\mu$ represents the mean solar irradiance, and $\sigma$ represents the standard deviation of solar irradiance.

The photovoltaic array generates current under sunlight through the photovoltaic effect, converting solar energy into electrical energy. The output power of the photovoltaic array is influenced by the solar irradiance $s$, and the following equation is used to convert solar radiation into solar energy:

$$P_{PV}(s) = \eta s A_s \tag{13}$$

Where $P_{PV}(s)$ represents the output power (kW) of the photovoltaic power generation system when the radiation level is $s$, $A_s$ is the surface area of the photovoltaic array (m$^2$), and $\eta$ is the power conversion efficiency of the photovoltaic system (PV).

### 2.2. Continuous abnormal weather classification based on scenario analysis

The probability distribution of continuous abnormal weather is an input parameter for the low-carbon park BIES optimization model. Compared to single abnormal weather events, the uncertainty of continuous abnormal weather is more significant. The coupling effects between different abnormal weather events enhance their temporal persistence, complexity, and impact range, posing greater challenges to the scheduling optimization of BIES. Therefore, it is necessary to decouple the coupling relationships of continuous abnormal weather to characterize the dependence characteristics between multidimensional random variables, introducing a dependency model based on Copula relationships [29]. This model primarily defines the correlation structure and marginal distributions between random variables, eliminates marginal effects through Copula, and then integrates them into a comprehensive joint distribution function.

Table 1. Classification of precipitation levels.

| Level | Hourly precipitation (mm) |
|---|---|
| Moderate rain | 0.4126~1.0416 |
| Heavy rain | 1.0417~2.0832 |
| Rainstorm | 2.0833~4.1666 |
| Torrential rain | 4.1666~10.4166 |

The uncertainty of continuous abnormal weather is represented by a scenario tree. This study primarily investigates scenarios of continuous torrential rain and strong wind, where wind speed is the random variable characterizing strong wind events, and precipitation is the random variable describing torrential rain events. According to the classification of precipitation levels and wind force grades by the China Meteorological Administration, precipitation is categorized into seven levels, while wind force is divided into seventeen levels. Under anomalous weather conditions, precipitation can be grouped into five levels, and wind force can be classified into ten levels. However, due to the excessive number of precipitation and wind speed combination scenarios, and the fact that the probability of some scenarios is close to zero, directly impacting the computational efficiency and accuracy of the BIES model, this study adopts a hierarchical optimization strategy. Specifically, scenarios with near-zero probability are disregarded, and based on actual anomalous weather data, the precipitation levels are simplified into four representative levels, while the wind force levels are reduced to four key levels. This approach significantly reduces computational complexity while ensuring the completeness of anomalous weather characteristics. The classification of precipitation and wind strength levels is shown in Table 1 and Table 2. As shown in Fig. 4, the combination of torrential rain and strong wind generates 16 different scenarios.

Table 2. Classification of wind strength.

| Wind Level | Name | Wind Speed (m/s) |
|---|---|---|
| 8 | Strong wind | 17.2~20.7 |
| 9 | Gale | 20.8~24.4 |
| 10 | Strong gale | 24.5~28.4 |
| 11 | Windstorm | 28.5~32.6 |

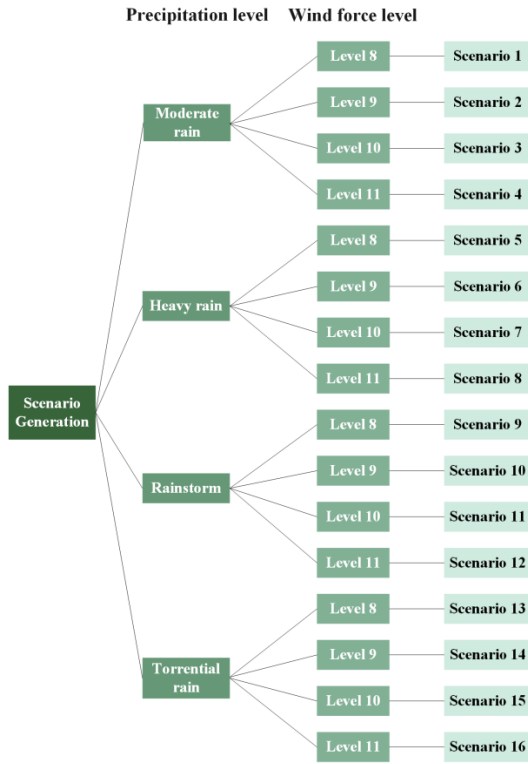

Fig. 4. Scenario tree structure for continuous abnormal weather based on stochastic programming.

## 2.3. Continuous abnormal weather correlation analysis

The dependency model based on Copula relationships can be used to quantify the uncertainty of the scale of continuous abnormal weather. Through collinearity, the dependency between multivariate distribution variables is probabilistically assessed, and by combining the occurrence probabilities of various abnormal weather events, the overall occurrence probability of continuous abnormal weather can be derived. Continuous abnormal weather characterized by variables such as wind speed and precipitation primarily includes strong winds and torrential rain. The wind speed model employs the Generalized Extreme Value (GEV) distribution, while the precipitation model uses a combination of the Generalized Pareto (GP) Distribution and Generalized Extreme Value distribution. The Generalized Pareto Distribution is expressed by Eq. (14), and the Generalized Extreme Value distribution is expressed by Eq. (15).

$$f_{(\mu,\sigma,\xi)}(x) = \frac{1}{\sigma}\left[1+\frac{\xi(x-\mu)}{\sigma}\right]^{\frac{1}{\sigma}-1} \tag{14}$$

$$f_{(\mu,\sigma,\xi)}(x) = \begin{cases} \exp\left(-\frac{x-\mu}{\sigma}\right)\exp\left[-\exp\left(-\frac{x-\mu}{\sigma}\right)\right] & \text{for } \xi = 0 \\ \left(1+\xi\frac{x-\mu}{\sigma}\right)^{-\left(1+\frac{1}{\xi}\right)}\exp\left[-\left(1+\xi\frac{x-\mu}{\sigma}\right)^{-\left(1+\frac{1}{\xi}\right)}\right] & \text{for } \xi \neq 0 \text{ and } \xi\frac{x-\mu}{\sigma} > 1 \\ 0 & \text{otherwise} \end{cases} \tag{15}$$

Where $\mu$ is the location parameter, $\sigma$ is the scale parameter, and $\xi$ is the shape parameter.

The wind speed and precipitation do not exhibit significant tail dependence characteristics, but they show strong correlation within the conventional numerical range. The optimal joint distribution between the two variables is represented by the Frank Copula. The Frank Copula function can describe the negative correlation between the variables. The distribution function and density function of the Frank Copula are as follows:

$$C_F(u,v) = -\frac{1}{\theta}\ln\left[1+\frac{(\exp^{-\theta u}-1)(\exp^{-\theta v}-1)}{\exp^{-\theta}-1}\right] \tag{16}$$

$$c_F(u,v) = -\frac{\theta(\exp^{-\theta}-1)\exp^{-\theta(u+v)}}{\left[(\exp^{-\theta}-1)+(\exp^{-\theta u}-1)(\exp^{-\theta v}-1)\right]^2} \tag{17}$$

## 2.4. Scenario generation and reduction

To effectively address the uncertainty in wind and photovoltaic output, a stochastic optimization method is employed, using scenario generation techniques to create a large number of wind and solar output scenarios that represent uncertain parameters. Scenario reduction techniques are then applied to approximate the accuracy of the original stochastic optimization with fewer scenarios, thereby reducing the computational complexity of the system.

Based on the aforementioned wind and solar power probability density distribution models, the Monte Carlo sampling method is used to generate a variety of typical scenarios that simulate the fluctuations in wind and solar output under different meteorological conditions. The generated scenarios not only cover various anomalous weather events but also capture the coupling relationships between these anomalous weather events and the probability weights of different scenarios.

By employing Monte Carlo simulation, random sampling is performed based on known mean and standard deviation to generate 2,000 scenarios, thereby simulating various uncertainty conditions. A random sampling method is applied to each scenario, and the results are stored in a matrix to characterize the system's performance under different conditions.

To reduce computational complexity and avoid the curse of dimensionality, a probability distance-based fast pre-generation elimination technique is applied. The retained scenarios must reflect the main variation trends of wind and solar output under anomalous weather conditions, ensuring the accuracy and efficiency of subsequent optimization.

By calculating the geometric distance between scenarios, it is possible to determine which scenarios have a high similarity. For each pair of scenarios, the sum of their absolute differences at each time point can be calculated to obtain the geometric distance:

$$d_{ij} = \sum_{t=1}^{24} |S_i(t) - S_j(t)| \quad (18)$$

Where $S_i(t)$ and $S_j(t)$ represent the values of the *i*-th and *j*-th scenarios at time *t*, respectively.

For each scenario, its average distance to all other scenarios is calculated as:

$$y_i = \frac{1}{m} \sum_{j=1}^{m} d_{ij} \quad (19)$$

Where $y_i$ represents the sum of the probability distances of the *i*-th scenario, and m is the total number of scenarios.

Using the scenario reduction method, the probability distribution of each scenario is first initialized to ensure equal initial probabilities. Subsequently, the average probabilistic distance between scenarios is computed, and the scenario with the smallest distance is selected as the representative scenario. This representative scenario is then merged with the nearest scenario, and the probability distribution is adjusted accordingly. The merged scenario is subsequently removed, and the probability vector and distance matrix are updated. This process is iteratively executed until the number of scenarios reaches the predefined threshold, thereby reducing computational complexity and improving optimization efficiency.

In the scenario generation problem of wind and solar power output, solar radiation and wind speed are considered the most critical weather uncertainty parameters. To accurately assess flexibility, this study employs a stochastic optimization method to model the aforementioned uncertainties. Additionally, to reduce computational complexity and shorten solution time, the unscented transform method is introduced to address the nonlinear function propagation issue in wind and solar output state estimation. The core idea of the unscented transform is to approximate the probability distribution of the original state by carefully selecting discrete sampling points, and then map these points to the new state space using a nonlinear function, thereby approximating the probability distribution of the target state.

In this study, solar radiation and wind speed, as two key uncertain parameters, are used as inputs for the original state distribution. By integrating the wind-solar power output uncertainty modeling function constructed in Section 2.1, a nonlinear mapping approach is employed to obtain the state distribution characteristics of wind and solar power output, thereby capturing the impact of uncertainty. Considering an uncertain nonlinear stochastic problem $y = f(x)$, where $y \in R^r$ is the uncertain output vector and $x \in R^q$ is the uncertain input vector, with a mean of $\mu$ and covariance $\sigma_x$. The symmetric elements of the matrix $\sigma_x$ represent the variances of the uncertain variables, while the asymmetric elements represent the covariances between different uncertain parameters. Using the unscented transform method, the mean $\mu_y$ and covariance $\sigma_y$ of the output variables are calculated as follows:

Step 1: Obtain $2q + 1$ samples from the input uncertain data:

$$x_0 = \mu \tag{20}$$

$$x_\omega = \mu + \left( \sqrt{\frac{q}{1-W^0}} \sigma_x \right), \omega = 1, 2, \ldots, q \tag{21}$$

$$x_\omega = \mu - \left( \sqrt{\frac{q}{1-W^0}} \sigma_x \right), \omega = 1, 2, \ldots, q \tag{22}$$

where $W_0$ is the weight associated with the mean value $\mu$.

Step 2: Evaluate the weight coefficients for each sample point:

$$W^0 = W^0 \tag{23}$$

$$W_\omega = \frac{1-W^0}{2q}; \quad \omega = q+1, \ldots, q \tag{24}$$

$$W_{\omega+q} = \frac{1-W^0}{2q}; \quad \omega+q = q+1, \ldots, 2q \tag{25}$$

$$\sum_{\omega \in s} W_\omega = 1 \tag{26}$$

Step 3: Sample $2q + 1$ points from the nonlinear function, and find the output samples using the following equation:

$$y_\omega = f(X_\omega) \tag{27}$$

Step 4: Evaluate the covariance $\sigma_y$ and mean $\mu_y$ of the output variable $y$.

$$\mu_y = \sum_\omega W_\omega y_\omega \tag{28}$$

$$\sigma_y = \sum_\omega W_\omega (y_\omega - \mu_y)(y_\omega - \mu_y)^T \tag{29}$$

## 3. Case Study

This study selects meteorological data from the European Centre for Medium-Range Weather Forecasts reanalysis dataset as the base input for generating the continuous abnormal weather dataset [30]. Using ten years of meteorological observation data, the wind-solar power output distribution under continuous anomalous weather conditions is constructed. To validate the proposed uncertainty modeling method in the context of BIES, typical anomalous weather days, such as heavy rainfall, strong winds, and windless conditions, are identified to support the wind-solar power output scenario generation in this study. Taking heavy rainfall as an example, meteorological data from Shenyang, Liaoning Province, China, from August 10 to 16, 2020, are selected as representative data for typical heavy rainfall abnormal weather. The geographical coordinates of Shenyang are 123°25′31″ and 41°48′11″, and its rainfall, wind speed, and surface horizontal radiation distributions are shown in Fig. 5, 6, and 7. For strong wind weather, meteorological data from Chifeng City, Inner Mongolia Autonomous Region, on March 16, 2024, are selected as representative data for typical strong wind abnormal weather. The geographical coordinates of Chifeng are 117°05′38″ and 43°29′18″. By analyzing these meteorological scenarios, this study evaluates the impact of continuous anomalous weather on wind-solar power output in BIES, thereby demonstrating the effectiveness of the proposed uncertainty modeling approach.

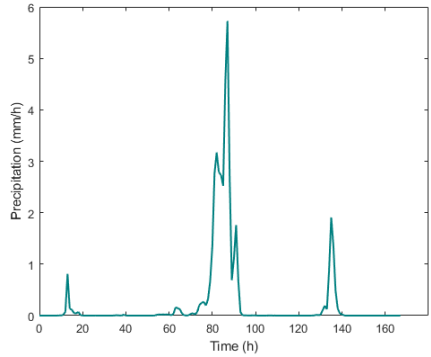

Fig. 5. Rainfall distribution on a typical day.

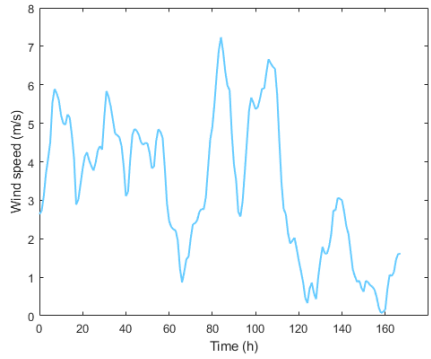

Fig. 6. Wind speed distribution on a typical day.

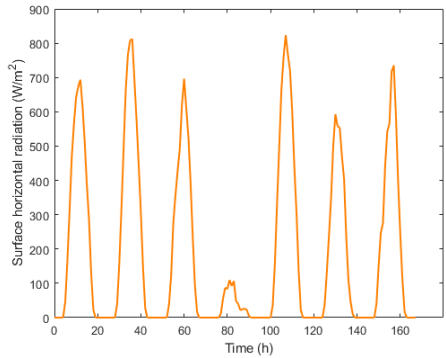

Fig. 7. Surface horizontal radiation distribution on a typical day.

In the simulation analysis, the wind turbine model Vestas-V90/2000 is selected, and the photovoltaic system adopts a fixed-panel lighting method. The simulation environment is set as an urban area. Table 3 lists the key parameters of the renewable energy equipment, including rated capacity and operation and maintenance costs.

Table 3. Equipment parameters of renewable energy.

| Type | Rated capacity $P_{i,max}$ /kW | Operation and maintenance costs $\lambda_{om,i}$ (CNY/kW) |
| --- | --- | --- |
| Photovoltaic system | 2000kW | 0.0235 CNY/kW |
| Wind turbine | 2000kW | 0.0196 CNY/kW |

## 4. Results and discussion

### 4.1. Generation of wind and photovoltaic output scenarios under normal weather conditions

Fig. 8 illustrates the reduced wind power and photovoltaic output scenarios under normal weather conditions. Initially, 2000 raw scenarios were generated based on the Monte Carlo simulation method, and then the number of scenarios was reduced to 5 using the fast predecessor elimination method based on probability distance. The results show that the reduced scenarios effectively capture the uncertainty characteristics of wind and photovoltaic outputs while ensuring the selection of typical scenarios. Specifically, the wind power output ranges from 0 to 1600 kW, with significant fluctuations primarily influenced by wind speed variations, demonstrating high randomness and instability. In contrast, the photovoltaic output fluctuates within the range of 0 to 2000 kW, driven by the periodic variations of solar radiation intensity, showing a more stable overall trend. The differences between the wind power output scenarios are more pronounced, indicating that the uncertainty in wind power is significantly higher than that in photovoltaic power, further validating the random disturbance effect of wind speed fluctuations on wind power output.

Fig. 9 presents the output characteristics of wind turbines and photovoltaic systems under normal weather conditions. From the output variation characteristics of a single scenario, the wind power output shows a multi-peak distribution over the course of the day, with a peak value of 1505.3 kW and a minimum output around 734.2 kW, indicating a significant influence of wind speed fluctuations on wind power output, leading to strong randomness. On the other hand, the photovoltaic output follows a typical single-peak trend, reaching its maximum between 12:00 and 14:00 and approaching zero in the early morning and evening, which aligns with the diurnal variation of solar radiation intensity. Overall, the uncertainty in wind power output is higher, posing a greater challenge for real-time power system dispatch. In contrast, photovoltaic output remains relatively stable, with its uncertainty mainly influenced by environmental factors such as cloud cover during short periods.

(a) (b)

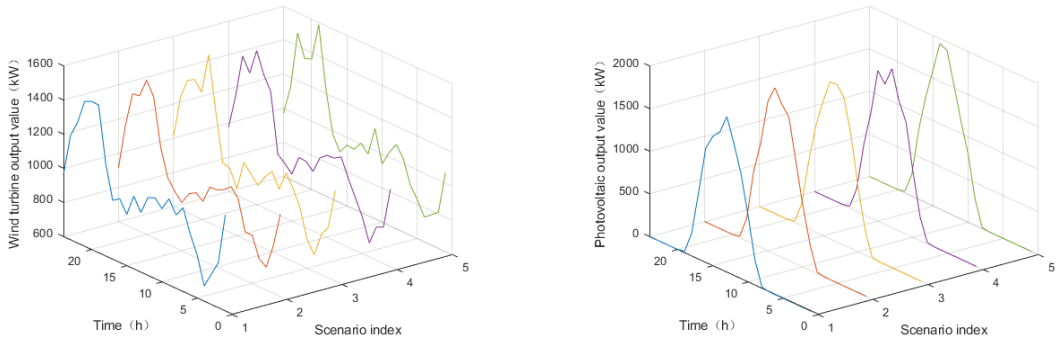

Fig. 8. Reduced wind turbine and photovoltaic output scenarios under normal weather conditions (a) wind turbine output (b) photovoltaic output.

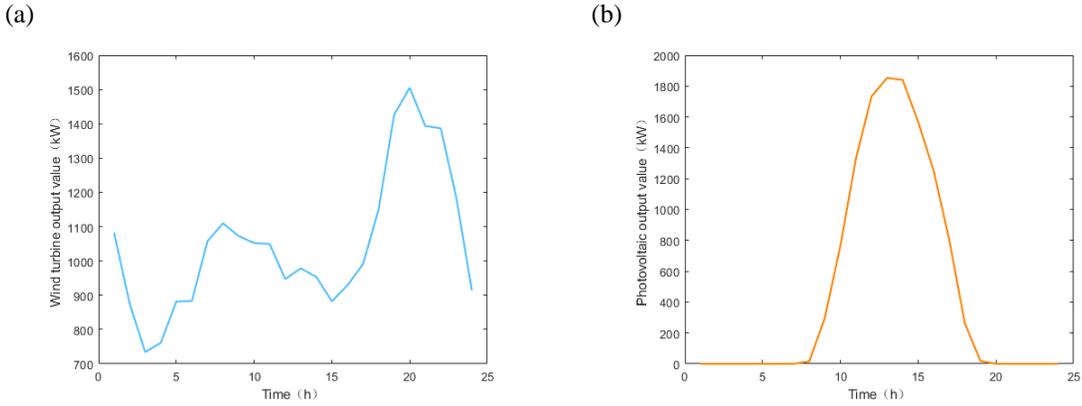

Fig. 9. Wind turbine and photovoltaic output under normal weather conditions (a) wind turbine output (b) photovoltaic output.

*4.2. Scenario generation of wind and solar power output on typical anomalous weather days*

Strong storms and heavy rainfall are typical examples of extreme weather, often accompanied by rapid variations in wind speed and significant fluctuations in solar radiation. These conditions directly impact the stability of wind and solar power output.

Under heavy rain conditions, photovoltaic power output fluctuates significantly, typically only 10% to 20% of the output under clear weather conditions. As long as there is light and visibility is not low, photovoltaic output will not be close to zero. The photovoltaic and wind power generation output under heavy rain conditions are shown in Fig. 10. Heavy rain is often accompanied by heavy precipitation, thunderstorms, and strong winds, which have a significant impact on wind power output. During the study period, windy weather mainly occurred between 20-114 hours and 155-167 hours, leading to increased wind power output.

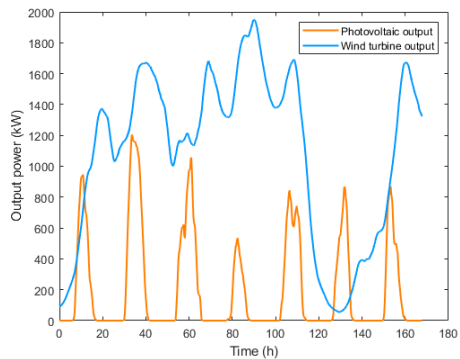

Fig. 10. Photovoltaic and wind turbine output under heavy rain conditions.

When wind speeds reach 17.2 m/s or higher, it is defined as strong wind conditions. In northeastern China, strong wind weather is more frequent in April, with an average of about 4 days of strong winds, and in some areas, it can reach 6-8 days. Strong wind conditions cause significant fluctuations in wind turbine output, while the impact on photovoltaic output is negligible. Therefore, this study mainly focuses on the wind turbine output characteristics under strong wind conditions. The wind turbine output curve for a typical strong wind day is shown in Fig. 11, where strong wind conditions persisted from 7-135 hours over six consecutive days. Due to variations in wind speed, wind turbine output fluctuates significantly. When the wind speed is lower, the output decreases accordingly, as shown by the lowest point in the figure. Since wind turbine output is generally high, excess output is stored in batteries for effective utilization.

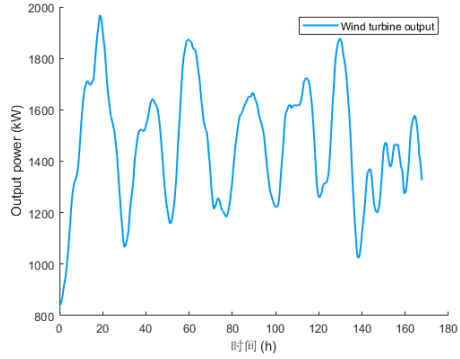

Fig. 11. Wind turbine output under strong wind conditions.

*4.3. Scenario generation of wind and solar power output in strong wind and rainstorm continuous abnormal weather*

When calculating the occurrence probability of consecutive anomalous weather events involving strong wind and heavy rain, the marginal distributions of each variable have a more significant impact on the results compared to the dependency functions between variables. The parameters of the marginal distribution functions for wind speed and precipitation are shown in Table 4. To characterize the correlation between wind speed and precipitation, this study employs the Frank Copula function for modeling and constructs the joint probability distribution surface in three-dimensional space along with its corresponding two-dimensional projection. Fig. 12 presents the joint probability density function of wind speed and precipitation. As observed in Fig. 12(b), strong wind and heavy rain are primarily concentrated in the lower-left quadrant, where the hourly precipitation is less than 5 mm and wind speed is below 30 m/s. A strong correlation is observed between the two variables, with wind speed exhibiting an increasing trend as precipitation intensifies.

Table 4. Parameters of marginal distribution functions for wind speed and precipitation.

| Variable | Distribution | Parameter | | |
|---|---|---|---|---|
| | | $\mu_1/\mu_2/\mu$ | $\sigma_1/\sigma_2/\sigma$ | $\xi_1/\xi_2/\xi$ |
| Wind speed | GEV | 11.892 | 8 | -0.175 |
| Precipitation | GP-GEV | 0.189 | 0.954 | 0.051 |
| | | 0.186 | 0.98 | 0.103 |

By constructing a scenario tree through the cross-combination of different types of wind speed and precipitation, the comprehensiveness of consecutive anomalous weather scenarios involving strong wind and heavy rain is enhanced. The occurrence probabilities of various scenarios are calculated using the Copula model, with specific values provided in Table 5. Fig. 13 illustrates the joint distribution samples generated by the Copula model, where the samples are distributed within the typical range of wind speed and precipitation values, exhibiting a relatively uniform distribution. Within this range, wind speed and precipitation demonstrate a certain degree of positive correlation. In the stochastic optimization process, each scenario is assigned a corresponding weight to ensure that the optimization model accurately reflects the wind-solar power output characteristics under different meteorological conditions. By incorporating historical anomalous weather data, the model is further aligned with real-world conditions, thereby improving the accuracy and robustness of the optimization results.

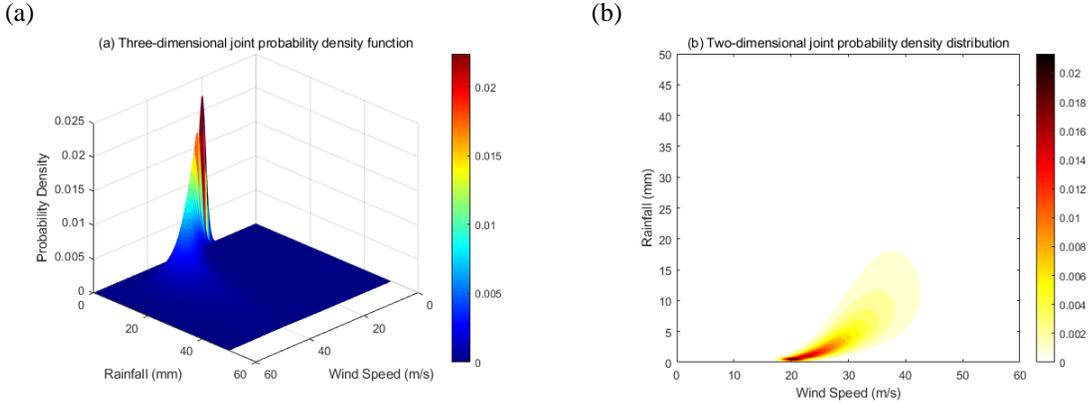

Fig. 12. Joint probability density function of wind speed and precipitation (a) three-dimensional space (b) two-dimensional plane.

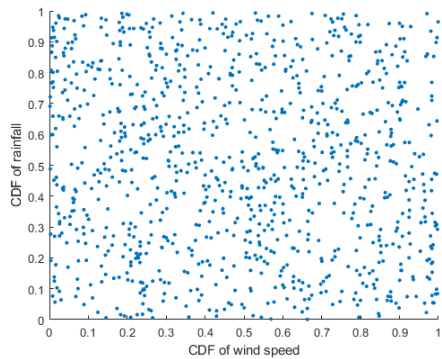

Fig. 13. Joint distribution samples generated by the Frank Copula.

The joint distribution exhibited by the Frank Copula demonstrates a strong tail dependence. Under extreme conditions, the correlation between wind speed and precipitation becomes more pronounced, with exceptionally high wind speeds often accompanied by significant precipitation. This characteristic makes the Frank Copula particularly suitable for describing wind speed and precipitation in anomalous weather events, especially in severe strong wind and heavy rain scenarios, as it effectively captures the strong tail dependence between these two variables.

Table 5. Probabilities of various scenarios in the scenario tree.

| Scenario | Precipitation level | Wind force level | Probability (%) |
| --- | --- | --- | --- |
| Scenario 1 | Moderate rain | Level 8 | 6.152 |
| Scenario 2 | Moderate rain | Level 9 | 8.682 |
| Scenario 3 | Moderate rain | Level 10 | 2.112 |
| Scenario 4 | Moderate rain | Level 11 | 0.251 |
| Scenario 5 | Heavy rain | Level 8 | 0.327 |
| Scenario 6 | Heavy rain | Level 9 | 5.315 |
| Scenario 7 | Heavy rain | Level 10 | 7.083 |
| Scenario 8 | Heavy rain | Level 11 | 1.684 |
| Scenario 9 | Rainstorm | Level 8 | 0.104 |
| Scenario 10 | Rainstorm | Level 9 | 3.053 |
| Scenario 11 | Rainstorm | Level 10 | 14.888 |
| Scenario 12 | Rainstorm | Level 11 | 11.802 |
| Scenario 13 | Torrential rain | Level 8 | 0.028 |
| Scenario 14 | Torrential rain | Level 9 | 0.965 |
| Scenario 15 | Torrential rain | Level 10 | 9.613 |
| Scenario 16 | Torrential rain | Level 11 | 27.941 |

Using the scenario generation and reduction method for consecutive anomalous weather of strong wind and heavy rain mentioned in Section 2.4, 16 representative scenarios were generated through Monte Carlo simulation and a rapid non-dominated elimination method based on probabilistic distance. The unscented transformation method was applied to handle the nonlinear function propagation in wind-solar power output state estimation. Fig. 14 illustrates the wind-solar power output distribution across the 16 representative scenarios. Two typical scenarios were selected from the generated scenario set, corresponding to the maximum and minimum wind-solar power output. Scenario 3 exhibits the highest total wind-solar power output, reaching 48,057.90 kW·h, which corresponds to strong wind with moderate rain conditions. Under strong wind conditions, wind speed remains within an optimal range without exceeding the turbine cut-out speed, allowing wind power generation to reach its maximum. Meanwhile, moderate rain has a relatively minor impact on solar irradiance, enabling photovoltaic generation to remain at a high level. In contrast, Scenario 13 represents the lowest total wind-solar power output at 36,800.30 kW·h, corresponding to strong wind with torrential rain conditions. Under these conditions, excessively high wind speeds exceed the turbine cut-out speed, triggering protective shutdown mechanisms and resulting in minimal wind power output. Additionally, torrential rain significantly reduces solar irradiance to its lowest level, further suppressing PV generation. These two scenarios can serve as wind and solar power output input data for the BIES bi-level collaborative optimization model to assess system optimization performance under various extreme conditions.

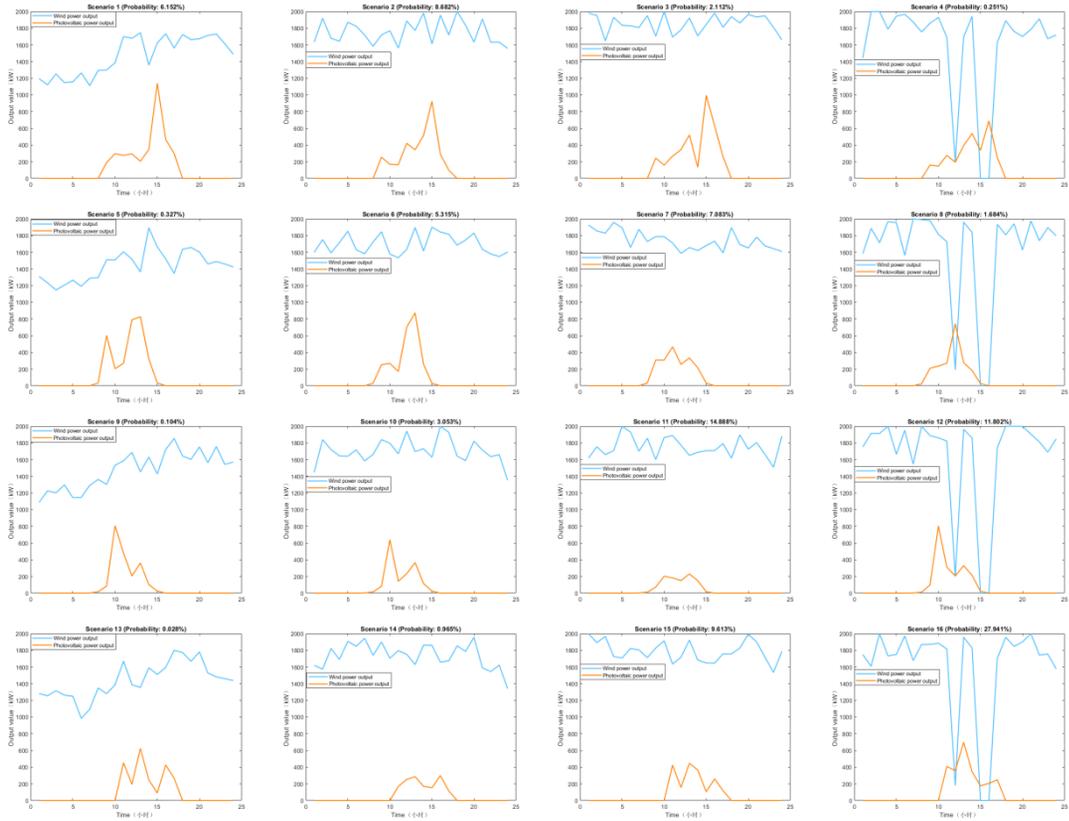

Fig. 14. Wind-solar power output across 16 representative scenarios.

*4.4. Evaluation metrics and validation methods*

*4.4.1. Marginal distribution function fitting test*

Precipitation and wind speed exhibit significant right skewness and heavy-tailed characteristics. The generalized extreme value distribution and the generalized Pareto distribution are commonly used for modeling heavy-tailed distributions. Therefore, the GEV and GP distributions were selected for fitting, and their goodness-of-fit was validated using the Akaike Information Criterion (AIC) test and the Kolmogorov-Smirnov (K-S) test.

Fig. 15 compares the fitting performance of the GP and GEV-GP models for precipitation distribution by analyzing the cumulative distribution functions of empirical and theoretical values. The GP model demonstrates a concentrated distribution fit, performing well in the 5-10 mm precipitation range but failing to capture the distribution of precipitation in the 0-5 mm range. In contrast, the GEV-GP model outperforms the GP model in fitting extreme precipitation values. The GEV-GP model retains the GP model's capability in describing tail behavior and extreme values while also achieving better fitting performance in the 0-10 mm range. Therefore, the GEV-GP distribution was selected for modeling precipitation distribution.

(a) 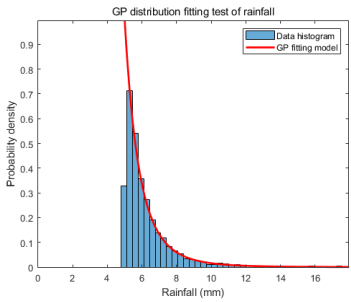
(b) 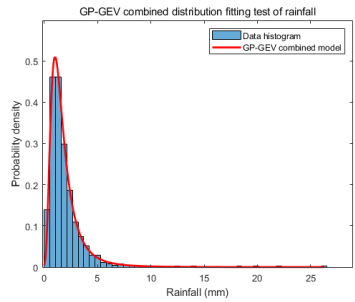

Fig. 15. Fitting performance validation of precipitation distribution using GP and GEV-GP models (a) GP distribution (b) GEV-GP distribution.

Fig. 16 presents the fitting performance of the marginal distributions for wind speed and precipitation. The left panel compares the cumulative distribution functions of empirical and theoretical values, while the right panel compares the probability density functions. The results indicate that the GEV distribution accurately captures the distribution characteristics of high wind speed values, demonstrating strong applicability under extreme wind speed conditions. Additionally, the GEV-GP distribution aligns more closely with the empirical curve in characterizing the tail behavior of precipitation, effectively capturing the distribution pattern of extreme precipitation events.

(a) 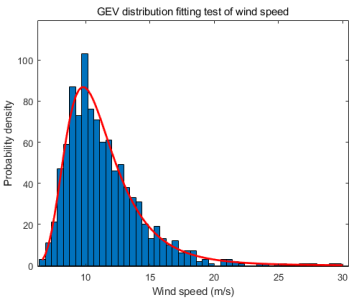
(b) 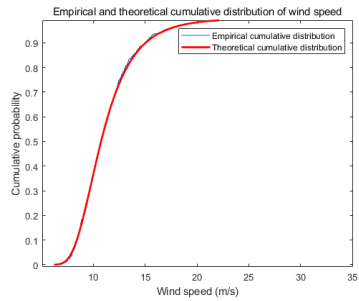
(c) 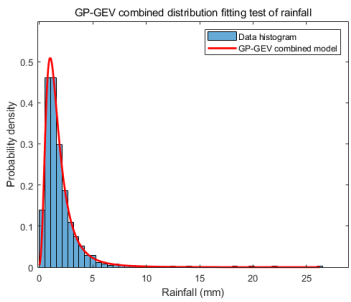
(d) 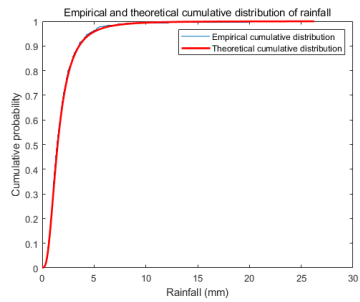

Fig. 16. Fitting performance validation of marginal distributions for wind speed and precipitation (a) cumulative distribution function of GEV distribution fitting (b) probability density function of GEV distribution fitting (c) cumulative distribution function of GEV-GP distribution fitting (d) probability density function of GEV-GP distribution fitting.

*4.4.2. Scenario reliability validation*

The reliability of the generated scenarios can also be assessed using quantile-quantile plots, which evaluate the similarity between two distributions. If the distribution of the generated scenarios closely matches that of the original scenarios, the quantile curve should align with a straight line. The Quantile-Quantile (QQ) plots for wind power output and photovoltaic power output are shown in Fig. 17 and Fig. 18, respectively. The quantile curves indicate that most data points are distributed along the diagonal, suggesting that the generated wind and PV power output scenarios largely preserve the statistical characteristics of the original scenarios. However, in the QQ plot of wind power output, certain deviations appear in the upper quantiles, where some data points diverge from the diagonal, implying that under extreme conditions, discrepancies may exist between the generated and original scenarios. This could be attributed to the influence of anomalous weather conditions, such as strong wind and heavy rain, which make it challenging for extreme wind power output values to fully align with the original scenario distribution. In contrast, the QQ plot of PV power output exhibits a higher degree of fit, with data points uniformly distributed around the diagonal, indicating that the rapid non-dominated elimination method based on probabilistic distance effectively preserves the statistical distribution characteristics of PV power output.

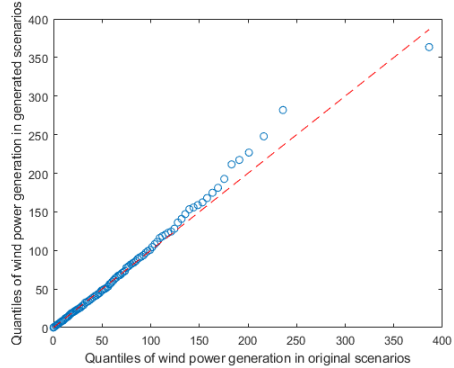

Fig. 17. Quantile-quantile plot of wind power output.

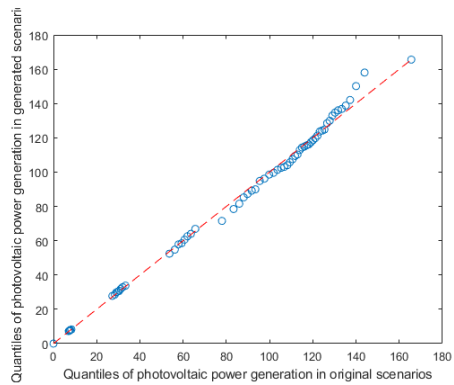

Fig. 18. Quantile-quantile plot of photovoltaic power output.

## 5. Conclusion

In stochastic optimization problems such as BIES optimization and dispatch, the accuracy of scenario generation plays a crucial role in system reliability assessment and the formulation of optimization strategies. Wind-solar power output is significantly influenced by weather conditions, and its uncertainty is further exacerbated under consecutive anomalous weather conditions. Consequently, the study of wind-solar power output scenario generation has attracted considerable attention. Properly modeling the coupling characteristics of wind speed and precipitation under consecutive anomalous weather conditions enables a more accurate characterization of wind-solar power output uncertainty, thereby improving the rationality of scenario generation and the effectiveness of system optimization. Building upon existing research on wind-solar power output uncertainty, this study further incorporates the coupling relationship between wind speed and precipitation under consecutive anomalous weather conditions and employs a scenario-based approach to model the uncertainty of wind-solar power output. First, a wind-solar power output model considering the joint distribution of wind speed and precipitation was established, with the Copula function used to describe their dependency structure. Then, Monte Carlo methods were applied to generate wind-solar power output scenarios incorporating consecutive anomalous weather characteristics, followed by scenario reduction based on probabilistic distance. The resulting wind-solar power output scenarios effectively capture the uncertainty of wind-solar power output under consecutive anomalous weather conditions while preserving the statistical properties of the original data.

It is worth noting that the Copula function demonstrates strong applicability in modeling the joint distribution of wind speed and precipitation. In this study, a static Copula model was employed to describe wind-solar power output uncertainty under consecutive anomalous weather conditions. However, the static Copula model assumes that the dependency structure remains stable over time, making it difficult to capture the dynamic variations in wind speed and precipitation under anomalous weather conditions. In particular, under extreme rainfall or strong wind conditions, non-stationarity and asymmetric dependencies between variables may lead to a decline in modeling accuracy. Future research will focus on dynamic modeling approaches such as time-varying Copula and hidden Markov Copula to enhance the model's adaptability to complex weather variations.

## CRediT authorship contribution statement

**Deyi Shao**: Writing – original draft, Methodology, Visualization, Validation, Software, Data curation, Conceptualization. **Hongru Li**: Writing – review & editing, Funding acquisition. **Jingsheng Li**: Writing – review & editing. **Xia Yu**: Writing – review & editing. **Xiaoyu Sun**: Writing – review & editing. **Bowen Han**: Writing – review & editing.

## Declaration of competing interest

The authors declare that they have no known competing financial interests or personal relationships that could have appeared to influence the work reported in this paper.

## Acknowledgements

The authors gratefully acknowledge the support by the National Natural Science Foundation of China (62473093).

## Data availability

The data presented in this study will be made available on request.